\def\ZZ         {{\mathbb Z}}

\def\CC         {{\mathbb C}}
\def\QQ         {{\mathbb Q}}

\def\ee         {{\rm e}}
\def\rk         {{\rm rk}}
\def\no         {\hspace{2pt}\colon}
\def\on         {\colon\hspace{2pt}}

\def\contr      {{\lrcorner}}
\def\Fock       {{\rm Fock}}

\def\dim        {{\rm dim}}

\documentclass[12pt]{amsart}
\usepackage{amsfonts}
\usepackage{amssymb}
\usepackage{amscd}
\usepackage{mathrsfs}

\newtheorem{lemma}{Lemma}[section]
\newtheorem{theorem}[lemma]{Theorem}
\newtheorem{corollary}[lemma]{Corollary}
\newtheorem{proposition}[lemma]{Proposition}
\theoremstyle{definition}
\newtheorem{definition}[lemma]{Definition}
\newtheorem{example}[lemma]{Example}
\newtheorem{question}[lemma]{Question}
\newtheorem{remark}[lemma]{Remark}
\theoremstyle{remark}
\newtheorem*{proof*}{Proof}

\begin{document}

\title{Chiral rings of vertex algebras of mirror symmetry}

\author{Lev A. Borisov}
\address{University of Wisconsin-Madison\\Mathematics Department\\ 480 Lincoln Dr, Madison, WI 53706, USA\\{\tt borisov@math.wisc.edu}.}

\begin{abstract}
{
We calculate chiral rings of the $N=2$ vertex algebras constructed
from the combinatorial data of toric mirror symmetry and show 
that they coincide with the description of stringy cohomology
conjectured previously in a joint work with A. Mavlyutov.
This constitutes an important reality check of the vertex 
algebra approach to mirror symmetry.
}
\end{abstract}

\maketitle

\section{Introduction}
Mirror symmetry was originally formulated by physicists 
as a correspondence between IIA models constructed for 
one family of Calabi-Yau varieties with IIB models constructed
for another family. IIA and IIB models are 
$N=(2,2)$ superconformal field
theories. While the very notion of $N=(2,2)$ superconformal field
theory is not yet properly axiomatized,
part of the structure of such a theory is a vertex algebra
with an $N=2$ superconformal structure, which does have a 
rigorous mathematical meaning. 
Every smooth Calabi-Yau variety $X$ gives rise to one 
such algebra, namely the cohomology of its chiral de Rham complex,
which is a certain sheaf of vertex algebras on $X$, constructed
in \cite{MSV}. 
Unfortunately, the definition
of chiral de Rham complex involves only the complex structure of 
the variety, not the metric or $B$-field. As such, it
lacks instanton corrections and, at best, is only related to the
large K\"ahler structure (radius) limit of the model.
In the case of 
the smooth Calabi-Yau variety, the chiral rings of the IIA
and IIB models are expected to give small quantum cohomology 
of the variety and the cohomology of the exterior powers of its tangent 
bundle. The chiral rings coming from the chiral de Rham
complex capture the latter ring, but only give the usual cohomology
product for the former.

Vertex algebras of mirror symmetry are a particular class
of vertex algebras with $N=2$ structure. Their definition in
\cite{borvert}
was 
inspired by the calculation of the cohomology of chiral de Rham
complex for the hypersurfaces in toric varieties, which is by
far the most common case of mirror symmetry. They are defined 
\emph{in purely combinatorial terms}
as a cohomology of certain explicit differentials on the 
lattice vertex algebras constructed from the combinatorial
data that underlie toric mirror symmetry. In
the hypersurface case, a given pair of dual
reflexive polytopes $(\Delta,\Delta^*)$
gives rise to a finite-dimensional
family of $N=2$ vertex algebras $V_{f,g}$ where $f$ and $g$ are
parameters that can be thought of as elements of finite-dimensional
vector spaces of complex-valued functions on the sets of lattice
points of $\Delta$ and $\Delta^*$ respectively. Parameter $f$
simply encodes the coefficients of the defining equation of
the hypersurface $X$ in the ambient
toric variety. Parameter $g$ 
can be thought of as a defining equation of an element $X^*$ of 
the mirror family of Calabi-Yau varieties, but its geometric 
meaning in terms of the original Calabi-Yau $X$ is far less clear.
In the smooth hypersurface case, after modding out by the symmetries 
induced by the torus action, the space of $g$ has the dimension
equal to $H^{1,1}(X)$. It is
\emph{conjecturally} related to the metric and $B$-field on $X$. 
In the more general case of $X$ given by a complete intersection
that corresponds to a nef-partition, the more appropriate language
is that of dual Gorenstein cones, see \cite{BatBor}. While 
the algebraic results of this paper are applicable to this more general
setting, in what follows we will focus on the hypersurface case
of toric mirror symmetry.

The basic premise of the vertex algebra approach to mirror
symmetry is that the family of vertex algebras $V_{f,g}$ 
is identical to the family of vertex algebras of IIA-IIB models
in the toric case, under some, yet unknown, identification of
the space of parameters $g$ and the complexified K\"ahler cone
of the Calabi-Yau hypersurface. At this moment this premise remains 
conjectural. Eventually, the hope is to give a mathematically rigorous
\emph{geometric} definition of the $N=2$ vertex algebras of 
IIA-IIB models for a given metric and $B$-field on an arbitrary
Calabi-Yau manifold and then calculate it in the toric hypersurface 
case. On the other hand, the family of vertex algebras $V_{f,g}$
is clearly compatible with the combinatorics
of mirror involution (i.e. a switch of $\Delta$ and $\Delta^*$
induces the mirror involution of the $N=2$ structure). It
has been also shown in \cite{borvert} that in the smooth ambient
toric variety case  for a given choice of parameters $f$
there is an algebraic family of vertex algebras
which has $V_{f,g}$ as its general fiber and has the cohomology 
of the chiral de Rham complex of $X$ as a special fiber. Here 
$g$ give the parameter space of the family. The author
hopes that the process of going from $V_{f,g}$ to 
the cohomology of chiral de Rham complex will eventually be 
interpreted as a mathematically rigorous counterpart of the 
physical procedure of going to the large radius limit.

This paper shows that vertex algebras $V_{f,g}$ pass 
two additional important reality checks. First, for a generic choice
of parameters $(f,g)$ the $N=2$ vertex algebra $V_{f,g}$
satisfies the positivity of Hamiltonians property expected 
of the vertex algebras of IIA-IIB models.

\smallskip\noindent
{\bf Theorem \ref{sigmatype}.}
\emph{For strongly non-degenerate $f$ and $g$ the $N=2$ vertex algebra 
$V_{f,g}$ is of $\sigma$-model type.}

Here we define an $N=2$ vertex algebra to be of $\sigma$-model type
if it satisfies the positivity, diagonalizability and 
finite-dimensionality of eigen-spaces 
of Hamiltonians $H_A$ and $H_B$, see Definition \ref{defsigma}.
We remark that in the case of smooth ambient toric variety 
this result can be established via the degeneration argument. 
Indeed, the statement
is obvious for the cohomology of the chiral de Rham complex
of the hypersurface $X_f$ given by a generic $f$.
Since the cohomology of the chiral de Rham complex fits into
an algebraic family with general fiber being $V_{f,g}$ 
(see \cite{borvert}), and
the dimensions of graded components in such family can only 
jump at special fibers, the result follows. The author thanks the
referee for the suggestion to stress this point. The more direct
approach of this paper has the advantage of avoiding the smoothness
condition which holds true only for a small minority of the 
reflexive polytopes, especially in higher dimensions.
In addition to giving a direct proof of Theorem
\ref{sigmatype}, our approach yields a more explicit description
of the chiral rings, which was the main motivation of the paper.

Our definition of chiral rings of $N=2$ vertex algebras of 
$\sigma$-model type follows closely that of 
\cite{LVW}. They are defined as zero eigenspaces of $H_A$ and
$H_B$ and are called $A$-ring and $B$-ring respectively. 
The main result of this paper is that for a general
choice of parameters $(f,g)$ the chiral rings of the algebra
$V_{f,g}$ have the expected graded dimension, as given by the 
(stringy) Hodge numbers of the hypersurface $X$. Moreover,
we give a rather explicit combinatorial description of the
fields in the lattice vertex algebra which descend to the 
fields in the chiral rings of $V_{f,g}$.

\medskip\noindent
{\bf Theorem \ref{main.chiral}.}
Let $f$ and $g$ be strongly non-degenerate.
Then both chiral rings of the 
$N=2$ vertex algebra 
$V_{f,g}$ are naturally isomorphic {\bf as vector spaces} to
the space
$$
W_{f,g}=\oplus_{\theta\subseteq K} R_1(\theta,f) \otimes_\CC
R_1(\theta^*, g).
$$

Here the direct sum is taken over the faces $\theta$
of the reflexive Gorenstein cone $K$ associated to $\Delta$,
$\theta^*$ denotes
the dual face of the dual cone and $R_1$ are vector spaces
defined by Batyrev in \cite{Bat.var} (see also \cite{anvar}). 
We remark that this 
paper implies that $W_{f,g}$  must possess two product structures, 
depending on whether we consider it to be the $A$-ring or
the $B$-ring of the algebra. Surprisingly,
we can not at present define these product structures 
in non-vertex-algebra terms.
On the other hand, the $A$-product structure on the ``diagonal
part'' which corresponds to $\theta = \{0\}$ is well-understood
and is illustrated in the famous quintic case, see
Example \ref{example.quintic}.

We remark that the spaces $W_{f,g}$ should be thought of as 
``(small) quantum stringy cohomology spaces'', under 
the $A$-ring product.
The first explicit occurrence of this idea in the literature
seems to be in \cite{MS}, which contains an easier
calculation of the chiral
rings of the related family of vertex algebras that correspond
to the situation when $X$ is a toric variety rather than a Calabi-Yau.
On the other hand, this idea is implicit in \cite{borvert}[Remark 8.6] 
and has been known to the author since then. 
Paper \cite{MS} also treats some Fano hypersurface examples.

We would like to clarify the relation of our construction 
with the concept of stringy Hodge numbers of Calabi-Yau 
hypersurfaces in toric varieties.
Toric mirror symmetry is best formulated for \emph{ample} 
hypersurfaces $X$ (i.e. without partial resolutions of singularities)
which are likely to be singular. There are two major reasons
for it. First, in dimension higher than three, Calabi-Yau 
hypersurfaces in toric varieties are unlikely to admit
a crepant resolution of singularities. Second, even in dimension
three, the crepant resolution of singularities is by no means
unique. In fact different crepant resolutions often produce
different cohomology rings. According to an idea of Batyrev,
in the singular case the usual cohomology
with the Hodge structure has to be replaced by some 
different double-graded vector space called stringy cohomology. 
While the graded dimension of the stringy cohomology vector space
has been defined in great generality (see 
\cite{Batyrev.cangor}), it is usually rather unclear
how to define the actual vector space. In fact, the author
has been suggesting for a while that it may not be possible 
to construct a single stringy cohomology vector space. Instead,
one should strive to
introduce a family of such vector spaces that depend on a finite
number of parameters.
It was observed in \cite{anvar} that in the case of hypersurfaces
in toric varieties, there is a natural construction of a
double-graded vector space of correct graded dimension, given
by an analog of $W_{f,g}$ above that takes into account the 
fan of the ambient toric variety.
While it was conjectured to be the mysterious stringy cohomology
space, it was not given a product structure expected
of such a space. The current paper rectifies this problem,
even though the product is not constructed explicitly.

While we use some key results of \cite{borvert} and \cite{anvar}, 
this paper is mostly self-contained. It is written almost entirely 
in algebraic terms, in hopes that it will be accessible to 
the specialists in the theory of vertex algebras. We are hoping
to make the case that vertex algebras of mirror symmetry are 
rich and beautiful objects that warrant further study. 

The paper is organized as follows.
In Section \ref{sec.cr} we recall the definitions of 
vertex algebras with $N=2$ structure and their chiral rings.
Section \ref{sec.lv} contains the definition of 
the particular class of lattice vertex algebras relevant to
this paper. In Section \ref{sec.vams} we recall the combinatorial
data of toric mirror symmetry and the definition of the
vertex algebras of mirror symmetry. Section \ref{sec.snd}
contains basic definitions and properties of strongly non-degenerate
coefficient functions which are a technical tool needed for
this paper. Our main results
are collected in Sections \ref{sec.main} and \ref{sec.sh}.
Finally, in Section \ref{sec.oq} we list some important
open questions. We hope that they will not stay open for long.

The author thanks the referee for suggestions on improving the
exposition. 

\section{N=2 vertex algebras and their chiral rings}\label{sec.cr}
In this section we recall the definitions of $N=2$ vertex algebras
and their chiral rings. We first define vertex algebras, according
to \cite{Kac}.

\begin{definition}\label{defvert}
A vertex algebra $V$ is first of all a super vector space over ${\CC}$,
that is  $V=V_0\oplus V_1$ where elements of $V_0$ are called bosonic or
even
and elements from $V_1$ are called fermionic or odd. In addition,
there is a fixed bosonic vector $\vert 0\rangle$ 
called {\it vacuum vector}.
The last part of the data that defines a vertex algebra is the
so-called state-field correspondence 
which is a parity preserving linear map from $V$ to ${\rm End}
V[[z,z^{-1}]]$
$$a\,\line(0,1){5}\hspace{-3.6pt}\to Y(a,z)=\sum_{n\in Z}a_{(n)}z^{-n-1}$$
such that for fixed $v$ and $a$ all $a_{(n)}v$ are zero for 
sufficiently big $n$.
This state-field correspondence must satisfy the following axioms.

\noindent$\bullet${\bf translation covariance:}
$\{T,Y(a,z)\}_-=\partial_z Y(a,z)$ where $\{,\}_-$ denotes
the usual commutator and $T$ is defined by $T(a)=a_{(-2)}\vert 0\rangle$;

\noindent$\bullet${\bf vacuum:} $Y(\vert 0\rangle,z)={\bf 1}_V,~
Y(a,z)\vert 0\rangle_{z=0}=a$;

\noindent$\bullet${\bf locality:} $(z-w)^N\{Y(a,z),Y(b,w)\}_{\mp}=0$
for all sufficiently big $N$, where $\mp$ is $+$ if and only if 
both $a$ and $b$ are fermionic. The equality is understood as an identity
of formal power series in $z$ and $w$.
\end{definition}

\begin{remark}
We will usually write $a(z)$ in place of $Y(a,z)$. The coefficients
$a_{(n)}$ are called \emph{modes} of the field $a(z)$.
For every two
elements $a$ and $b$ there is  an operator product expansion (OPE)
$$
a(z)b(w)=\sum_{i=1}^N \frac{c^i(w)}{(z-w)^i}+
\no a(z)b(w)\no$$
where the meaning of the symbols $\frac 1 {(z-w)^i}$ 
and $\no\no$ in the above formulas is as in Chapter 2 of \cite{Kac}. 
Operator product expansions is a convenient way to encode the 
(super-)commutators of the modes of $a$ and $b$.
\end{remark}

\begin{definition}
An $N=2$ vertex algebra is a vertex algebra $V$ with the following
additional structure. There are fixed bosonic elements $L$ and $J$
and fixed fermionic elements $G^+$ and $G^-$ of $V$ such that 

\noindent$\bullet$ $L_{(0)}=T$

\noindent$\bullet$ $L_{(1)}$ is diagonalizable and satisfies
$\{L_{(1)},Y(a,z)\}_-=z\partial_zY(a,z)+Y(L_{(1)}a,z)$ for all $a\in V$.

\noindent$\bullet$ Fields $L(z), J(z), G^\pm(z)$ satisfy the OPEs
$$
\begin{array}{rcl}
L(z)L(w)&=& \frac{c/2}{(z-w)^4} + \frac{2L(w)}{(z-w)^2}
+\frac{\partial_wL(w)}{z-w} + \no L(z)L(w) \on,
\\
L(z)J(w)&=&\frac{J(w)}{(z-w)^2} + \frac {\partial_w J(w)}{z-w} + 
\no L(z)J(w)\on,
\\
L(z)G^{\pm}(w)&=&\frac{(3/2)G^{\pm}(w)}{(z-w)^2} + \frac{\partial_w G^{\pm}
(w)}{z-w}+ \no L(z)G^{\pm}(w)\on,
\\
J(z)J(w)&=&\frac{c/3}{(z-w)^2}+\no J(z) J(w)\on,
\\
J(z)G^{\pm}(w)&=&\pm\frac{G^{\pm}(w)}{z-w}+ \no J(z)G^{\pm}(w)\on,
\\
G^{\pm}(z)G^{\mp}(w)&=&\frac{2c/3}{(z-w)^3}\pm\frac
{2J(w)}{(z-w)^2}
+\frac{2L(w)\pm\partial_wJ(w)}{z-w}+ \no G^{\pm}(z)G^{\mp}(w)\on ,
\\
G^{\pm}(z)G^{\pm}(w)&=&\no G^{\pm}(z)G^{\pm}(w)\on.
\end{array}
$$
Here $c$ is a constant and we will call $\hat c=c/3$ 
the \emph{central charge} of the $N=2$ vertex algebra.
\end{definition}

Let $V$ be an $N=2$ vertex algebra. We are particularly 
interested in the operators $L_{(1)}$ and $J_{(0)}$. OPE
of $L$ and $J$ implies that $L_{(1)}$ and $J_{(0)}$ 
commute with each other. In all the situations that we
will consider in this paper they will provide $V$ with a double grading.
The following definition is inspired by \cite{LVW}.

\begin{definition}\label{defsigma}
We call an $N=2$ vertex algebra $V$ an \emph{$N=2$ vertex algebra
of $\sigma$-model type} if the eigenspaces $V_{\alpha}$ of $L_{(1)}$ 
are finite-dimensional for all $\alpha$ and 
zero except for $\alpha\in\frac12\ZZ_{\geq 0}$,
and the operators $H_A:=L_{(1)}-\frac 12 J_{(0)}$ and 
$H_B:=L_{(1)}+\frac 12 J_{(0)}$ have only nonnegative
integer eigenvalues.
\end{definition}

\begin{remark}
From the physicists' point of view the conditions above are
satisfied for $\bar\partial$-quotients of the $N=(2,2)$ algebras 
constructed (empirically) from sigma-models on Calabi-Yau manifolds. 
The central charge of such algebra equals the dimension
of the Calabi-Yau manifold.
The operators $H_A$ and $H_B$ are the Hamiltonians for the two topological
twists of the $\sigma$-model, which explains the notations. 
In the case of smooth Calabi-Yau manifold, the $A$-ring should equal
its small quantum cohomology while the $B$-ring should equal the
cohomology of exterior powers of its tangent bundle.
\end{remark}

We now define chiral $A$-rings and chiral 
$B$-rings of $N=2$ vertex algebras 
of $\sigma$-model type as zero eigenspaces of $H_A$ and $H_B$ 
respectively. Letters $A$ and $B$ here stand for 
the rings of the topological quantum field theories of IIA
and IIB models.
The ring structures are provided by the following proposition.
We formulate it for the $A$-ring with $B$-ring case being completely
analogous.

\begin{proposition}
Let $V$ be an $N=2$ vertex algebra of $\sigma$-model type.
Let $a,b\in V$ be such that $H_A(a)=H_A(b)=0$.
Then all the modes of $a$ and $b$ super-commute and 
$c(z):=a(z)b(z)$ is a field of $V$ which corresponds 
to an element $c$ with $H_A(c)=0$. This provides the kernel of $H_A$
with the structure of an associative super-commutative ring.
\end{proposition}

\begin{proof}
The proof is essentially contained in \cite{LVW}.
Consider the OPE 
$$
a(z)b(w)=\sum_{i=1}^N \frac{c^i(w)}{(z-w)^i}+
\no a(z)b(w)\no
$$
of $a$ and $b$.
A standard calculation shows that the $H_A$ grading of 
$c^i$ is negative, hence $c^i=0$ for all $i$. As a result,
the modes of $a(z)$ and $b(w)$ super-commute. The field 
$a(z)b(z)$ corresponds to the element $a_{(-1)}b_{(-1)}\vert 0\rangle$ 
which is easily shown to have $H_A=0$, since $a_{(-1)}$ and
$b_{(-1)}$ preserve the $H_A$ grading.
\end{proof}

\begin{remark}
We recall that mirror involution is the automorphism of $N=2$
vertex algebra that does not change the vertex algebra $V$
but moves $L, J, G^+, G^-$ to $L, -J, G^-, G^+$. Indeed, this 
change preserves the OPEs of $N=2$ algebra. This involution
switches $A$-rings and $B$-rings.
\end{remark}

\begin{remark}
From physical considerations, one often expects to find a vector 
space isomorphism between the $A$-ring and the $B$-ring
of a vertex algebra of $\sigma$-model type,
given by the so-called \emph{spectral flow} (see \cite{LVW}).
We will establish this isomorphism explicitly
for the $N=2$ vertex algebras considered in
this paper, see Remark \ref{remspectral}. 
It is important to keep in
mind that spectral flow \emph{does not preserve
the ring structure.}
\end{remark}

\section{Lattice vertex algebras}\label{sec.lv}
In this section we will recall the definition of lattice vertex algebras
in a particular case. Let $M$ and $N$ be two dual free abelian groups
which will be fixed throughout this section. We denote 
the natural pairing of $m\in M$ and $n\in N$ by $m\cdot n$ and extend
it by linearity to a pairing between $M_\CC$ and $N_\CC$.

The infinite
Heisenberg algebra associated to $M\oplus N$ is the associative
$\CC$-algebra with generators $m^{bos}_i$ and  $n^{bos}_j$ where 
$m\in M_\CC$, $n\in N_{\CC}$, $i,j\in\ZZ$. 
The relations are linear relations on $m^{bos}_i$ for a fixed $i$
and $n^{bos}_j$ for fixed $j$ according to the relations in $M_\CC$
and $N_\CC$, as well as the relations
$$
\{m^{bos}_i,n^{bos}_j\}_-=i(m\cdot n)\delta_{i+j}^0,~
\{m^{bos}_i,m^{bos}_j\}_-=\{n^{bos}_i,n^{bos}_j\}_-=0.
$$
Here $\{\}_-$ denotes the commutator and $\delta$ is the Kronecker symbol.

Analogously, the infinite Clifford algebra associated to $M\oplus N$ has 
generators $m^{ferm}_i$ and  $n^{ferm}_j$ for $m\in M_\CC$,
$n\in N_\CC$ and $i,j\in \ZZ+\frac 12$. In addition to linear relations 
for fixed $i$ or $j$, there are anti-commutator relations
$$
\{m^{ferm}_i,n^{ferm}_j\}_+=(m\cdot n)\delta_{i+j}^0,
\{m^{ferm}_i,m^{ferm}_j\}_+=\{n^{ferm}_i,n^{ferm}_j\}_+=0\hskip -1.6pt.
$$

\begin{remark} The superscripts $bos$ and $ferm$ stand for bosonic 
and fermionic respectively.
\end{remark}

We will consider the following representation of the direct sum 
of the Clifford and Heisenberg algebras. Let 
${\rm Fock}_{{\bf 0}\oplus{\bf 0}}$ be the tensor product over all 
integers $i<0$ of the symmetric algebras of $M_\CC$ and  $N_\CC$
tensored by the product over all half-integers $i<0$ of the exterior
algebras of $M_\CC$ and $N_\CC$. 
For $i<0$ the action of the generators $m^{bos}_i$ and $n^{bos}_i$ 
on ${\rm Fock}_{{\bf 0}\oplus{\bf 0}}$
of the infinite Heisenberg algebra is defined as the multiplication 
by $m^{bos}_i$ and $n^{bos}_i$ respectively. For $i>0$ the action of
$m^{bos}_i$ is defined as $i$ times the corresponding differentiation
on the symmetric algebra of $N_\CC$ for the index $-i$, and similarly
for $n^{bos}_i$. The action of the generators of infinite Clifford
algebra is similarly defined by either multiplication or differentiation.
Finally, the action of the generators $m^{bos}_0$ and $n^{bos}_0$
(which commute with all other generators and with each other) is defined 
to be zero. As will be apparent later, this is the reason for the subscript
${\bf 0}\oplus {\bf 0}$ in our notations for the space. Traditionally,
the constant $1$ in this product of symmetric and exterior powers is
denoted $\vert {\bf 0},{\bf 0}\rangle$.

It is well-known that the space 
${\rm Fock}_{{\bf 0}\oplus {\bf 0}}$ has a structure of 
a vertex algebra with $\vert {\bf 0},{\bf 0}\rangle$
as its vacuum vector and the parity given by 
the number of the modes of $m^{ferm}$ and $n^{ferm}$ modulo $2$.
The fields of this algebra are 
linear combinations of normal ordered products of the 
fields 
$$
\begin{array}{l}
m^{bos}(z):=\sum_{i\in\ZZ} m^{bos}_iz^{-i-1}\\
n^{bos}(z):=\sum_{i\in\ZZ} n^{bos}_iz^{-i-1}\\
m^{ferm}(z):=\sum_{i\in\ZZ+\frac 12} m^{ferm}_iz^{-i-\frac 12}\\
n^{ferm}(z):=\sum_{i\in\ZZ+\frac 12} n^{ferm}_iz^{-i-\frac 12}
\end{array}
$$
and their derivatives with respect to $z$. Normal ordered product means 
that the generators with positive modes are always applied first
(see \cite{Kac}). We remark that, in particular, the fields $m^{bos}(z)$ 
correspond to the elements $m^{bos}_{-1}\vert {\bf 0},{\bf 0}\rangle$
which we will call $m^{bos}$, abusing the notations. Similarly, 
we define $n^{bos}$, $m^{ferm}$ and $n^{ferm}$ as the elements 
that correspond to the above fields. We remark that while we 
have $m^{bos}_i=m^{bos}_{(i)}$ and $n^{bos}_i=n^{bos}_{(i)}$
where the brackets in the subscripts signify the convention from
Definition \ref{defvert}, we get $m^{ferm}_i=m^{ferm}_{(i-\frac12)}$
and $n^{ferm}_i=n^{ferm}_{(i-\frac12)}$.

We have not yet used the lattice structure of $M$ and $N$, since 
we have only worked with their complexifications. We will now  
define a bigger vertex algebra ${\rm Fock}_{M\oplus N}$ which 
will contain ${\rm Fock}_{{\bf 0}\oplus{\bf 0}}$ as a subalgebra.
As a vector space, ${\rm Fock}_{M\oplus N}$ is a tensor product over
$\CC$ of ${\rm Fock}_{{\bf 0}\oplus{\bf 0}}$ and the group algebra
$\CC[M\oplus N]$. We will denote the tensor product of $\vert
{\bf 0},{\bf 0}\rangle$ with the basis element of the group algebra
that corresponds to $m\oplus n\in M\oplus N$ by $\vert m,n\rangle$.
We will denote the corresponding subspace of ${\rm Fock}_{M\oplus N}$
by ${\rm Fock}_{m \oplus n}$, in agreement with our convention for
${\rm Fock}_{{\bf 0}\oplus {\bf 0}}$.
Every element of ${\rm Fock}_{m\oplus n}$ can be 
thought of as a polynomial in the generators of the infinite Heisenberg
and Clifford algebras with negative modes applied to $\vert m,n\rangle$.

\begin{theorem}
The vector space ${\rm Fock}_{M\oplus N}$ has a natural structure
of a vertex algebra.
\end{theorem}

We omit the proof, since it is standard. However we will
describe the fields of this algebra below to fix the notations.

We extend the action of $m^{bos}(z)$ and $n^{bos}(z)$ to
${\rm Fock}_{M\oplus N}$ by requiring that $m^{bos}_0$ and 
$n^{bos}_0$ act by $(m\cdot n'){\bf id}$ and 
$(m'\cdot n){\bf id}$ on ${\rm Fock}_{m'\oplus n'}$ for all 
$m'\in M,n'\in N$. The fermionic fields $m^{ferm}$ and $n^{ferm}$
extend without any changes. As a result, the action of all
fields of ${\rm Fock}_{{\bf 0}\oplus{\bf 0}}$ extends to
${\rm Fock}_{M\oplus N}$. 
To describe the vertex algebra structure on ${\rm Fock}_{M\oplus N}$
we also need the 
so-called \emph{vertex operators} which will be the fields of 
the vertex algebra that correspond to the elements $\vert m,n\rangle$.
\begin{definition}(\cite{Kac})
Let $\gamma_{m+n}$ denote the multiplication by $[m\oplus n]$
in the semigroup algebra $\CC[M\oplus N]$. Then the field
$\no\ee^{\int (m^{bos}+n^{bos})(z)}\on$ is defined by 
$$
\gamma_{m+n}
(-1)^{m^{bos}_0}
z^{m^{bos}_0+n^{bos}_0}
\prod_{i<0}\ee^{-(m^{bos}_i+n^{bos}_i)\frac {z^{-i}}{i}}
\prod_{i>0}\ee^{-(m^{bos}_i+n^{bos}_i)\frac {z^{-i}}{i}}
$$
and is  called \emph{vertex operator} for $m\oplus n$.
Here $z^{m^{bos}_0}$ acts on ${\rm Fock}_{m'\oplus n'}$
by $z^{m\cdot n_0}$ and similarly for $z^{n^{bos}_0}$ and 
$(-1)^{m^{bos}_0}$.
\end{definition}

A general field of the vertex algebra ${\rm Fock}_{M\oplus N}$
is a linear combination of normal ordered products 
of one the fields $\ee^{\int (m^{bos}+n^{bos})(z)}$ and
the fields of ${\rm Fock}_{{\bf 0}\oplus{\bf 0}}$ where the 
normal ordering is taken with respect to the generators of 
${\rm Fock}_{{\bf 0}\oplus{\bf 0}}$.

\begin{remark}
We will often omit normal ordering from our formulas to 
simplify the notations, since it will be present in all calculations.
\end{remark}

\begin{definition}
For any subset $S$ of $M\oplus N$ we define the subspace
${\rm Fock}_{S}$ of ${\rm Fock}_{M\oplus N}$
as $\oplus_{m\oplus n\in S}{\rm Fock}_{m\oplus n}$.
If $S$ is a subsemigroup, then ${\rm Fock}_{S}$ is 
a vertex subalgebra.
\end{definition}

\section{Vertex algebras of mirror symmetry}\label{sec.vams}
In this section we recall the combinatorics of the toric mirror
symmetry and define certain families of $N=2$ vertex algebras in 
terms of this combinatorics. These algebras have been first
introduced in \cite{borvert}.

The fundamental notion of the toric mirror symmetry is that of
a pair of dual reflexive Gorenstein cones. As in the previous
section, we fix a pair of dual lattices $M$ and $N$. We recall that
a rational polyhedral cone is the intersection of 
a finite number of closed rational halfspaces. We will also assume
that a cone does not contain any sublattices. The dual of a cone $K$
is defined as the set of all elements of the dual rational vector space 
that are nonnegative on $K$.
\begin{definition}(\cite{BatBor})
Let $K$ and $K^*$ be rational polyhedral cones in $M$ and $N$
respectively that are dual to each other. These cones are called
reflexive Gorenstein cones iff there exist lattice elements $\deg\in K$
and $\deg^*\in K^*$ such that the lattice 
generators $m$ of all one-dimensional
faces of $K$ satisfy $m\cdot \deg^*=1$ and the 
lattice generators $n$ of all
one-dimensional faces of $K^*$ satisfy $\deg \cdot n=1$. The number
$\deg\cdot \deg^*$ is called \emph{index} of the pair of reflexive 
Gorenstein cones.
\end{definition}

\begin{definition}
Let $K$ and $K^*$ be dual reflexive Gorenstein cones. We denote 
by $\Delta$ (resp. $\Delta^*$) the set of points of $m\in K\cap M$ 
(resp. $K^*\cap N$) that satisfy $m\cdot \deg^*=1$ (resp. $\deg\cdot n=1)$.
These are polytopes of dimension $\rk M-1$.
\end{definition}

\begin{remark}
For the readers familiar with the mirror symmetry for hypersurfaces 
in toric varieties, we remark that in the case of index $1$ the 
polytopes $\Delta$ and $\Delta^*$ are reflexive (see \cite{Batyrev}).
In a more general case of a Calabi-Yau complete intersection of $k$ 
hypersurfaces, associated to a nef-partition, 
one gets reflexive Gorenstein cones of index $k$ (see \cite{BatBor}).
Although not all reflexive Gorenstein cones come from this construction,
their number is known to be finite in any given dimension.
\end{remark}

To define vertex algebras of mirror symmetry we will need additional
data, namely a pair of \emph{coefficient functions} $f\colon\Delta\to \CC$
and $g\colon\Delta^*\to \CC$ which assign a constant to every lattice
point of $\Delta$ and $\Delta^*$. The terminology is justified by 
the fact that in the hypersurface case these functions give the
coefficients of the defining equation.

\begin{definition}
A coefficient function $f$  is called \emph{non-degenerate} iff 
the quotient of the semigroup ring $\CC[K]$ by the ideal generated
by the elements $z_n=\sum_{m\in\Delta}f(m)(m\cdot n)[m]$ for $n\in N$
is finite-dimensional. Here $[m]$ denotes the element of $\CC[K]$
that corresponds to the element $m\in K$.
Non-degenerate functions $g$ are defined 
similarly. 
\end{definition}

One can show (see \cite{Bat.var} or \cite{locstring})
that the non-degenerate
functions $f$ form a Zariski open subset in the set of all possible $f$.
Moreover, for any basis $\{n_k\}$ of $N_\CC$ and any non-degenerate $f$,
elements $z_{n_k}$ form a regular sequence in $\CC[K]$. 

We are now ready to describe the vertex algebras of interest.

\begin{definition}
Let $K$ and $K^*$ be dual reflexive Gorenstein cones. Let
$f$ and $g$ be non-degenerate coefficient functions for $K$ and $K^*$
respectively. 
Then the \emph{vertex algebra of mirror symmetry} 
$V_{f,g}$ is the cohomology 
of the lattice vertex algebra ${\rm Fock}_{M\oplus N}$ with respect to
the operator $D_{f,g}=D_f+D_g$ where 
$$
D_f:={\rm Res}_{z=0} \sum_{m\in\Delta}f(m)m^{ferm}(z)\ee^{\int m^{bos}(z)},
$$
$$
D_g:={\rm Res}_{z=0} \sum_{n\in\Delta^*}g(n)n^{ferm}(z)\ee^{\int n^{bos}(z)}.
$$
\end{definition}

We remark that an easy OPE calculation shows that 
$D_f^2=D_g^2=D_fD_g+D_gD_f=0$, so the above operator is indeed 
a differential. 
The vertex algebra structure on ${\rm Fock}_{M\oplus N}$ 
induces a vertex algebra structure on the cohomology by $D_{f,g}$.
Moreover, one can endow $V_{f,g}$ with $N=2$ structure by considering
the elements $G^+, G^-, J, L$ that correspond to the fields
\begin{equation}\label{n2str}
\begin{array}{rcl}
G^+(z)&=&\sum_k (n^k)^{bos}(z)(m^k)^{ferm}(z)-\partial_z\deg^{ferm}(z)
\\
G^-(z)&=&\sum_k (m^k)^{bos}(z)(n^k)^{ferm}(z)-\partial_z(\deg^*)^{ferm}(z)
\\
J(z)&=&\sum_k (m^k)^{ferm}(z)(n^k)^{ferm}(z)+\deg^{bos}(z)-
(\deg^*)^{bos}(z)
\\
L(z)&=&\sum_{k}(m^k)^{bos}(z)(n^k)^{bos}(z) 
+\frac12\sum_k\partial_z(m^k)^{ferm}(z)(n^k)^{ferm}(z)\\
&&
\hskip -17.2pt 
-\frac12\sum_k
(m^k)^{ferm}(z)\partial_z(n^k)^{ferm}(z)
-\frac12\deg^{bos}(z)-\frac12(\deg^*)^{bos}(z)
\end{array}
\end{equation}
where $\{m^k\}$ and $\{n^k\}$ are (any) dual bases of $M_\CC$ 
and $N_\CC$. 
\begin{proposition}\label{ntwo}
Elements $G^\pm, J,L$ lie in the kernel of $D_{f,g}$ and therefore
descend to its cohomology. The central charge of the
resulting $N=2$ algebra is $\rk M-2\deg\cdot\deg^*$.
\end{proposition}

\begin{proof}
This is a standard calculation in OPEs and is left to the reader.
\end{proof}

\section{Strongly non-degenerate coefficient functions}\label{sec.snd}
The setup for this section is the following.
Let $K$ be a cone in $M$ such that there exists an 
element $\deg^*\in N={\rm Hom}(M,\ZZ)$ such that the
values of the linear function $\deg^*:M\to \ZZ$ 
on the minimum lattice generators
of one-dimensional cones of $K$ are $1$. We will
denote the natural pairing $M\times N\to \ZZ$ by $\cdot$.
In the applications
of this paper the cone $K$ will be reflexive Gorenstein,
but we only need the above weaker condition. 
We denote by $\Delta$ the set of lattice points  $m\in K$
that satisfy $m\cdot \deg^*=1$.
As before, $\deg^*$ provides the semigroup ring $\CC[K]$
with a grading and every coefficient function $f:\Delta\to\CC$ 
gives a homogeneous element of degree one in $\CC[K]$.

We recall that 
for any basis $\{n_i\}$ of $N_\CC$ and any non-degenerate 
$f$ the elements $z_{n_i}$ form a regular sequence. 
We now introduce a technical notion of strongly
non-degenerate coefficient functions, which will be used in
the next section. Let $n_0$ be an element of $N$. It provides
$\CC[K]$ with additional grading $\cdot n_0$
which may or may not be nonnegative. Let $\{n_i\},i=1,\ldots,\rk N$
be a basis of $N$. We define $z_{i,j}$ to be the $j$-th graded
component of the element $z_i:=z_{n_i}$ with respect to the 
$\cdot n_0$ grading. We consider all values of $\cdot n_0$ that 
occur for some $m\in\Delta$, but some of the $z_{i,j}$ may well
be zero. We will denote the set of possible values 
of $m\cdot n_0$ by $I(n_0)$.

We then consider the Koszul complex 
associated to the sequence of elements $z_{i,j}\in \CC[K]$
with $j\in I(n_0)$.
It is given by 
\begin{equation}\label{Koszul}
\wedge (\oplus_i\oplus_{j\in I(n_0)}
\CC e_{i,j})
\otimes_\CC \CC[K]
\end{equation}
where $\wedge$ means the exterior algebra and the differential
is given by $d=\sum_{i,j} \contr(e_{i,j}) z_{i,j}$.
Clearly, the cohomology $W(f;n_0)$ of this complex is independent from 
the choice of the basis $\{n_i\}$. 

\begin{proposition}
For every non-degenerate $f$ and every $n_0$ the space 
$W(f;n_0)$ is finite-dimensional. 
\end{proposition}

\begin{proof}
Denote by  $j_{min}$ the 
smallest element in $I(n_0)$. 
The vector space 
$E=\oplus_i\oplus_{j\in I(n_0)}\CC e_{i,j}$
is the direct sum of the subspace $E_1$ generated 
by $e_i = e_{i,j_{min}}$ 
and the subspace $E_2$ generated
by $e_{i,j}'=e_{i,j}-e_{i,j_{min}}, j>j_{min}$.
The direct sum decomposition
$E=E_1\oplus E_2$ induces the decomposition
$$
\wedge^r E = \oplus_{k+l=r} \wedge^k E_1 \otimes_\CC \wedge^l E_2.
$$
We use this decomposition to realize the Koszul complex 
associated to the sequence of elements $z_{i,j}\in \CC[K]$
as the total complex of a double complex. Namely, we can decompose
the differential $d$ on $\wedge E\otimes_\CC\CC[K]$
as $d=d_1+d_2$ where 
$$
d_1 = \sum_{i} \contr(e_{i}) z_i,~d_2 = \sum_{i,j>j_{min}(n_0)}
\contr(e_{i,j}') z_{i,j}.
$$
Here we have used $z_i = \sum_j z_{i,j}$. 

In view of the spectral sequence for the double complex,
to show that $W(f;n_0)$ is finite-dimensional, it is enough
to show that the cohomology of $\wedge E_1\otimes_\CC\wedge
E_2\otimes_\CC\CC[K]$ with respect to $d_1$ is finite-dimensional.
This cohomology is isomorphic to the tensor
product over $\CC$ of $\wedge E_2$ and the cohomology 
of  $\wedge E_1\otimes_\CC\CC[K]$ with respect to $d_1$. The 
latter is simply the Koszul complex for $\{z_i\}$, so the 
statement follows from the assumption that 
$f$ is non-degenerate.
\end{proof}

\begin{definition}
A non-degenerate $f$ is called \emph{
strongly non-degenerate} if for each $n_0$ the dimension of $W(f,n_0)$
is minimum among all possible choices of the coefficient functions.
\end{definition}

\begin{proposition}
Strongly non-degenerate functions form a Zariski open subset of 
all coefficient functions.
\end{proposition}

\begin{proof}
Indeed, while there are infinitely many different choices of
$n_0$ they lead to only finitely many choices for $\{z_{i,j}\}$,
since $\Delta$ has only finitely many points. For each such choice
the minimality of the dimension of the quotient is a Zariski open
condition.
\end{proof}

Observe that we have not used the double grading by 
$\cdot \deg^*$ and $\cdot n_0$ in the 
definition of strongly non-degenerate $f$. However, $W(f,n_0)$
does inherit a double grading from $\CC[K]$ since the Koszul complex
is taken for a sequence of homogeneous elements. There is a 
unique, up to a shift, way to endow the Koszul complex with the 
structure of the double-graded module over $\CC[K]$, in a way to make 
its differential grading-preserving. 
\begin{proposition}\label{dg}
Let $f$ be strongly non-degenerate. Then for every $n_0$ the dimension
of each double-graded piece of $W(f,n_0)$ is the smallest among all 
possible choices of $f$.
\end{proposition}

\begin{proof}
It is clear that the dimension of each graded component 
$W(f,n_0)$ is at least as big as it is for a general choice of $f$.
However, the sum of the dimensions is the same as for the general 
coefficient function, which means that all graded dimensions are
minimum.
\end{proof}

We will also need to consider \emph{partial}
semigroup rings\footnote{These have been called deformed
semigroup rings in our earlier papers.}. 
Namely, let ${\mathcal T}$ be a regular triangulation of $\Delta$.
The regularity condition is best described in terms of the 
corresponding decomposition $\Phi$ of $K$ into simplicial cones.
It means that there is a real valued continuous function $h$ on $K$ which
is linear on each cone and satisfies the convexity relation
$$h(x+y)\geq h(x)+h(y)$$
for all $x,y\in K$ where the equality holds iff $x$ and $y$ lie
in the same cone of $\Phi$. Moreover, $h$ can be picked to have integer
values on the lattice points of $K$.
We can define the \emph{partial
semigroup ring} $\CC[K]^\Phi$ by redefining the product of 
monomials $[m][m_1]$ to be $[m+m_1]$ if $m$ and $m_1$ are in
the same cone of $\Phi$ and zero otherwise. It has been shown
in \cite{anvar} that $\CC[K]^\Phi$ is a Cohen-Macaulay ring.
Moreover, if we introduce the notion of non-degenerate coefficient
function for $\CC[K]^\Phi$, such functions form a Zariski open
subset  (see \cite{anvar}). As a consequence, the theory developed
in this section transfers immediately to $\CC[K]^\Phi$ 
in place of $\CC[K]$. We will write $^\Phi$ to indicate the 
partial semigroup theory.

We will be interested in the relation between the dimensions
of the double-graded components of $W(f^\Phi,n_0)^\Phi$ and $W(f,n_0)$. 
\begin{proposition}\label{deform}
Let $f$ be a strongly non-degenerate coefficient
function for $\CC[K]$ and let $f^\Phi$ be any non-degenerate 
coefficient function for $\CC[K]^\Phi$.
Then 
for a fixed $n_0$ and a fixed value of the double grading, the dimension
of the double-graded component of $W(f,n_0)$ does not 
exceed that of the corresponding component of $W(f^\Phi,n_0)^\Phi$.
\end{proposition}

\begin{proof}
The regularity of triangulation allows us to see $\CC[K]^\Phi$
as a certain limit of $\CC[K]$. 
Namely, consider a family of multiplication structures $*_q$ on 
the vector space $\CC[K]$ by defining 
$$[m]*_q[m_1]=[m+m_1]q^{h(m+m_1)-h(m)-h(m_1)}$$
where $h$ is the function from the definition of regularity.
For $q\neq 0$ the resulting ring is isomorphic to the ring $\CC[K]$
under the rescaling $[m]\to[m] q^{h(m)}$. At $q=0$ the resulting
ring is precisely $\CC[K]^\Phi$.

Consider the family of Koszul complexes with $f_q=f^\Phi$. Since 
$f^\Phi$ is non-degenerate for $\CC[K]^\Phi$, the dimension
of the cohomology at $q=0$ is finite. Since this construction
is compatible with the double grading which has finite-dimensional
components, it implies that the dimensions of the graded components
of the cohomology of the Koszul differential at $q\neq 0$ do not 
exceed those for $q=0$. On the other hand, the Koszul complex 
at $q\neq 0$ is isomorphic to the Koszul complex for $\CC[K]$
and the coefficient
function $f_q(m)=f^\Phi(m)q^{-h(m)}$. Finally, by Proposition \ref{dg},
the dimensions of the graded pieces of the 
cohomology of the Koszul complex for $f$ do not exceed those
for any other coefficient function, so in particular
they do not exceed those for $f_q$.
\end{proof}

\section{Main Theorem}\label{sec.main}
The goal of this section is to show that for general choices of 
$f$ and $g$ the vertex algebra $V_{f,g}$ is of $\sigma$-model type.
We will also describe the subspaces of the fields of the lattice
vertex algebra which descend to the chiral rings of $V_{f,g}$,
although we defer the more explicit calculation of $A$- and $B$-rings
of the algebra until the next section.

The following proposition describes the action on
${\rm Fock}_{M\oplus N}$ of the operators 
$L_{(1)}$ and $J_{(0)}$ of the $N=2$ structure defined in 
Section \ref{sec.vams}.
\begin{proposition}\label{lj}
Consider an element 
$$v=\prod_p (m^p)^{bos}_{-i_p}\prod_q (n^q)^{bos}_{-j_q}
\prod_r (m^r)^{ferm}_{-i'_r}\prod_s (n^s)^{ferm}_{-j'_s}
\vert \hat m,\hat n
\rangle$$
where $m^p,m^r\in M_\CC$ and $n^q,n^s\in N_\CC$.
Then 
$$
\begin{array}{rcl}
L_{(1)}v &=& (\hat m\cdot \hat n+\frac12\deg\cdot \hat n+\frac12\hat m
\cdot \deg^*\\
&&
+\sum_pi_p+\sum_qj_q+\sum_ri'_r+\sum_sj'_s)v
\\
J_{(0)}v &=& (\deg\cdot \hat n -\hat m\cdot \deg^*+\sum_r1-\sum_s1)v.
\end{array}
$$
\end{proposition}

\begin{proof}
The zeroeth modes of $\deg^{bos}(z)$ and $(\deg^*)^{bos}(z)$ act
on ${\rm Fock}_{\hat m\oplus \hat n}$ by multiplication by 
$\deg\cdot\hat n $ and $\hat m\cdot \deg^*$ respectively, which
accounts for the terms with $\deg$ and $\deg^*$ in the above formulas.
The rest is a standard calculation in free bosonic and fermionic 
fields (see \cite[Sections  3.5, 3.6]{Kac}).
\end{proof}

\begin{corollary}\label{hb}
In the notations of the above proposition, the operator $H_A$ 
is given by 
$$
H_Av = (\hat m\cdot(\deg^*+\hat n)
+\sum_pi_p+\sum_qj_q+\sum_r(i'_r-\frac12)+\sum_s(j'_s+\frac12))v.
$$
\end{corollary} 

\begin{proof}
Operator $H_A$ was defined as $L_{(1)}-\frac 12 J_{(0)}$, so it remains
to use Proposition \ref{lj}.
\end{proof}

We first give some estimates of the cohomology of ${\rm Fock}_{K\oplus N}$
by $D_f$. While $H_A$ has some negative eigen-values on 
${\rm Fock}_{K\oplus N}$ due to the $\hat m\cdot(\deg^*+\hat n)$
term in Corollary \ref{hb}, it turns out to be nonnegative on
the $D_f$ cohomology, at least for a strongly non-degenerate $f$.

\begin{proposition}\label{dfprop} For a strongly non-degenerate $f$
the $D_f$-cohomology of ${\rm Fock}_{K\oplus N}$ has 
only eigenspaces of $H_A\geq 0$. Moreover, the $H_A=0$ eigenspace 
comes from
${\rm Fock}_{K\oplus (K^*-\deg^*)}$.
\end{proposition}

\begin{proof}
Clearly, $D_f$ does not change the $N$-grading. 
Let us fix 
for a moment $n_0\in N$ and work in $V_{n_0}={\rm Fock}_{K\oplus n_0}$.
If $n_0$ lies in $K^*-\deg^*$ then $H_A$ is nonnegative on 
$V_{n_0}$. As a result, it is enough to show that for $n_0\not\in
K^*-\deg^*$ the values of $H_A$ on $D_f$ cohomology of $V_{n_0}$
are at least $1$.

From now on we will assume that $n_0\not\in K^*-\deg^*$.
The space $V_{n_0}$ is filtered by the spaces $W_k$
defined as the span of all elements $v$ of Proposition \ref{lj}
with $\sum_p1-\sum_q1\geq k$, i.e. the number of bosonic modes
of coming from lattices $M$ minus the number of bosonic modes
coming from lattice $N$ is at least $k$.
Since $D_f$ is a linear combination of products of some modes of 
$m^{ferm}$ and some modes of $\ee^{\int m^{bos}}$, it preserves 
the filtration.
Together with $\cdot\deg^*$ grading, this provides $V_{n_0}$
with the structure of the filtered complex, with $D_f$ as its 
differential.

Notice that the filtration by $W_k$ is also compatible with $H_A$ grading,
so it is enough to consider a fixed eigenvalue of $H_A$.
For every value of $\cdot\deg^*$, there are only finitely many possible 
$m\in K$. As a result, Corollary \ref{hb} shows that the dimension
of the corresponding eigenspace of $H_A$ in the fixed $\cdot\deg^*$
graded component of ${\rm Fock}_{m\oplus n_0}$ is finite.
As a result, the $\{W_k\}$ filtration is finite, which
assures the convergence of the spectral sequence of the filtered complex.
So to show that a certain $H_A$-component of $D_f$ cohomology of 
$V_{n_0}$ is zero, it is enough to show that the cohomology of
the differential $d$ induced by $D_f$ on the $H_A$-component of 
$W_k/W_{k+1}$ is zero for all $k$. 

Let us consider the action of $d$ on $\oplus_kW_k/W_{k+1}$. Recall that 
$$
D_f=\oint dz \sum_{m\in \Delta}f(m) m^{ferm}(z)\ee^{\int m^{bos}(z)}.
$$
The only mode of $\ee^{\int m^{bos}(z)}$ that acts non-trivially on
$\oplus_kW_k/W_{k+1}$ is given by the shift $\gamma_m$. 
As a result, the action of $d$ on $\oplus_kW_k/W_{k+1}$ is given by 
$$d=\sum_{m\in\Delta} f(m) m^{ferm}_{l_m}\gamma_m$$
where $l_m$ depends on 
$m$ and $n_0$ only. In fact, 
\begin{equation}\label{lm}
l_m=m\cdot(n_0 +\deg^*)-\frac 12=m\cdot n_0+\frac 12.
\end{equation}
It can be derived, for instance, from the fact that $D_f$ commutes 
with $H_A$ and therefore does not change the $H_A$ grading.

As before, we denote by $I(n_0)$ the set of all possible values
of $m\cdot n_0$ for $m\in\Delta$. 
The action of $d$ on $\oplus_kW_k/W_{k+1}$
really comes from its action on the space
$$U_{n_0}=
\Big(\wedge(\oplus_{j\in I(n_0)}F_j)\Big)\otimes_\CC \CC[K].
$$ 
Here
$F_j$ is $(M_\CC)^{ferm}_{j+\frac 12}$ for $j+\frac 12< 0$ and 
$(N_\CC)^{ferm}_{-j-\frac 12}$ otherwise. 
Namely, the space $\oplus_kW_k/W_{k+1}$
is the tensor product of $U_{n_0}$ and the polynomial algebra in
all other non-positive modes of 
$m^{bos},n^{bos},m^{ferm},n^{ferm}$, and the action of
$d$ comes from the $U_{n_0}$ component of this tensor product.
Hence, the cohomology of $d$ is the tensor product of its cohomology
on $U_{n_0}$ and this other space. Clearly, extra non-positive modes
can not decrease the value of $H_A$ so it is enough to
show that the cohomology of $d$ on $U_{n_0}$ has the desired
properties with respect to $H_A$ grading. At this point one can
forget about vertex algebras and talk simply about the action of
$d=\sum_{m\in\Delta} f(m) [m]\otimes m^{ferm}_{l_m}$ on the space 
$U_{n_0}$.
The action of $m^{ferm}_{j+\frac 12}$ 
is either multiplication or contraction,
depending on whether or not $j+\frac 12 < 0$.

The key observation is that this complex is precisely isomorphic to the 
Koszul complex considered in Section \ref{sec.snd}.
To see this,
we consider the isomorphism between $\otimes_{j\in I(n_0)} \wedge F_j$ 
and $\wedge (\oplus_i\oplus_{j\in I(n_0)}
\CC e_{i,j})$ of \eqref{Koszul} given as follows.
Let $\{n_i\}$ be a basis of $N$ and let $\{m_i\}$ be the dual
basis of $M$. For $j+\frac 12>0$
we can identify $F_j$ with $\oplus_i\CC e_{i,j}$ 
by mapping $(n_i)^{ferm}_{-j-\frac 12}$ to $e_{i,j}$.
Then the part of the differential $d$ 
\begin{equation}\label{partofd}
\sum_{m\in\Delta,l_m=j+\frac 12} f(m) [m]\otimes m^{ferm}_{l_m}
\end{equation}
corresponds exactly to the part of the Koszul differential
\begin{equation}\label{partofK}
\sum_{i} \contr e_{i,j}z_{i,j}.
\end{equation}
For $j+\frac 12 < 0$ we identify $\wedge F_j\cong 
\wedge (M_\CC)^{ferm}_{j+\frac 12}$ with 
$\wedge (\oplus_i \CC e_{i,j})$ by mapping 
$$(m_{i_0})^{ferm}_{j+\frac 12}\wedge \ldots 
\wedge (m_{i_s})^{ferm}_{j+\frac 12}$$
to
$$
(\contr e_{i_0,j})\ldots (\contr e_{i_s,j})
\Big(e_{1,j}\wedge\ldots\wedge e_{\rk N,j}\Big).$$
This has an effect of changing the multiplication action
of $m^{ferm}_{j+\frac 12}$ into the contraction action
by a linear combination of $e_{i,j}$, and once again
\eqref{partofd} translates into \eqref{partofK}.
We will now use Proposition \ref{deform}.

Recall that $n_0$ does not lie in $K^*-\deg^*$, which means that
at least one of the $l_m$ is negative. Let us fix this $m$ and call
it $\hat m$. Consider the regular triangulation
of $\Delta$ so that $\hat m$ is contained in any nonboundary simplex.
Such a triangulation can be constructed from the following collection
of heights on lattice points of $\Delta$. Assign $\hat m$ height $0$ and
assign all other lattice points of $\Delta$ heights that are generic
and close to $1$. Clearly, every simplex of maximum dimension of
the resulting triangulation of $\Delta$ contains $\hat m$.
Moreover, every nonboundary simplex of this triangulation contains 
$\hat m$. Otherwise, 
it can be extended to a nonboundary simplex of codimension at most 
one with the same property. Such a codimension one simplex is a 
boundary of two simplices of maximum dimension, only one of which
can contain $\hat m$. We denote by $\Phi$ the fan on $K$ that corresponds
to this triangulation. This triangulation allows us to redefine 
$U_{n_0}$ by redefining the product structure on $\CC[K]$
as in Section \ref{sec.snd}.
We will denote the corresponding space by $U_{n_0}^\Phi$. The new
differential will still be denoted by $d$, abusing the notations slightly. 
By Proposition \ref{deform}, it is enough to show the 
vanishing of $H_A\leq 0$ cohomology
of $U_{n_0}^\Phi$. In fact, it is enough to show this for one specific
value of $f$, namely the one that has values $1$ for all vertices of 
$\Delta$ and values $0$ otherwise.

We denote by $\Phi(k)$ the set of nonboundary 
cones of $\Phi$ of dimension $k$.
We can write $U_{n_0}^\Phi$ as the cohomology of the complex $\mathcal E$
$$
0\to E^{\rk M}\to E^{\rk M-1}\to \ldots\to 0
$$
where 
$$E^{k}
=\oplus_{C\in \Phi(k)}{\CC[C]\otimes\prod_{l_m} \wedge \CC^{\rk M}}$$
and the differential $d'$ is constructed from the restriction maps
with signs coming from some choices of orientation on the cones.
Indeed, for every $m\in K$ the part of this complex with this
$M$-grading comes from nonboundary cones of $\Phi$ that contain $m$.
It is isomorphic to a complex that calculates reduced homology of 
a sphere times $\prod_{l_m}\wedge\CC^{\rk M}$. More importantly,
the action of $d$ on the cohomology of 
$U_{n_0}^\Phi$ is induced from the
following action on the above complex. The action of $d$ on 
$\CC[C]\otimes \prod_{l_m}\wedge\CC^{\rk M}$ is
\begin{equation}\label{donC}
\sum_{i=1}^{\dim C} [m]\otimes m_{i,l_i}
\end{equation}
where $m_i$ are the generators of one-dimensional faces of $C$,
$l_i$ is a shorthand for $l_{m_i}$ 
and $m_{i,l_i}$ act by either multiplication or by contraction.
It is easy to see that $d$ and the differential $d'$ of $\mathcal E$
can be combined to get a double complex. 

Since the cohomology of $\mathcal E$ with respect to $d'$ 
is concentrated at the first column, its
cohomology with respect to the total differential $d+d'$ is the
same as the cohomology of $U_{n_0}^\Phi$ with respect to $d$.
On the other hand, there is a spectral sequence that converges from
the cohomology of the total space of $\mathcal E$ with respect to
$d$ to the $d$-cohomology of $U_{n_0}^\Phi$. As a result, it is
sufficient to check that the $d$-cohomology of $\CC[C]\otimes 
\prod_{l_m}\wedge \CC^{\rk M}$ by the differential of \eqref{donC}
has zero $H_A$ eigenspaces for non-positive eigenvalues.
If we recall how the $H_A$ grading is defined, it is easy to see that
it is sufficient to consider the case of $\dim C=\rk M$. Indeed, 
in the general case we can see that the complex is a product of a
complex that involves the exterior algebras of $\CC$-span of $C$ and its
dual and a vector space of nonnegative $H_A$ grading. Analogously, 
we can ignore all $\wedge \CC^{\rk M}$ for $l_m$ that are not among
$l_i$. 

Let's denote the space $\CC[C]\otimes \prod_{l_i}\wedge\CC^{\rk M}$
by $V$. Denote by $\{n_i\}$ the basis of $N_\QQ$ dual to $\{m_i\}$.
As a graded space with a differential, $V$ is isomorphic
to $W\otimes \CC[C]\times \prod_{i}\wedge\CC e_i$ where $W$ is some
non-negatively graded vector space with the differential coming from $d$
and $e_i$ stands for either $m_i$ or $n_i$ depending on the sign of
$l_i$. This complex is isomorphic to $W$ times the Koszul complex 
for the regular sequence $\{[m_i]\}$. 
As a result, the $d$-cohomology is given by $W$ tensored 
with $\oplus_{m\in {\rm Box}(C)}\CC[m]$, which is 
the quotient of $\CC[C]$ by the ideal generated by $\{[m_i]\}$.
${\rm Box(C)}$ is defined as the set of all $m\in M$ that
have coordinates in $[0,1)$ in the basis $\{m_i\}$.
However, there is an additional shift in the grading due to the 
fact that isomorphism with the Koszul complex occurs only after 
switching from the multiplication to the contraction. This 
adds $\sum_{l_i<0} (-l_i-\frac12)$ to the eigenvalues of $H_A$. 

If $m=\sum_i a_im_i$, then the value of $H_A$ on the element $[m]$ is 
\begin{equation}\label{technical}
\sum_ia_im_i\cdot (n_0+\deg^*)-\sum_{l_i<0} (l_i+\frac 12)
\geq 
\sum_{l_i<0} (a_i-1)(l_i+\frac12)
\end{equation}
where we have used $m_i\cdot (n_0+\deg^*)=l_i+\frac 12$.
This number is positive since by our construction there is at least one
$l_i\leq -\frac 32$, namely the one that corresponds to $\hat m$.
\end{proof}

\begin{remark}
A slightly weaker version of Proposition \ref{dfprop} is proved
in Section 9 of \cite{borvert}. However, the degeneration argument
presented there only shows this statement for a generic rather than
general $f$, i.e. one may have to exclude countably many 
Zariski closed subvarieties. The argument of this paper is more direct
and assumes only that $f$ is strongly non-degenerate. We don't
know if the statement holds for arbitrary non-degenerate $f$, although
we have no examples to the contrary.
\end{remark}

We can generalize Proposition \ref{dfprop} to other eigen-values of 
$H_A$.
\begin{proposition}\label{somebigr}
Let $f$ be strongly non-degenerate. Then for every $k>0$ the
corresponding eigenspace of $H_A$ on $D_f$-cohomology of 
${\rm Fock}_{K\oplus N}$ comes from 
${\rm Fock}_{K\oplus (K^*-r\deg^*)}$ for some sufficiently big $r$.
\end{proposition}

\begin{proof}
We follow the proof of Proposition \ref{dfprop}. At the very 
last step we observe that if $r$ is big enough, then $l_i$ that 
corresponds to $\hat m$ will be sufficiently negative to assure that 
$\sum_{l_i<0} (a_i-1)(l_i+\frac12)$ in
\eqref{technical} is bigger than $k$.
\end{proof}

\begin{theorem}\label{main}
Let $f$ and $g$ be strongly non-degenerate.
Then $D_{f,g}$-cohomology of ${\rm Fock}_{M\oplus N}$
has only nonnegative integer eigenvalues of $H_A$. 
Moreover, the $H_A=0$ eigenspace comes from 
${\rm Fock}_{K\oplus (K^*-\deg^*)}$. The operator $H_B$
also has only nonnegative integer eigenvalues on $V_{f,g}$
and its kernel comes from ${\rm Fock}_{(K-\deg)\oplus K^*}$.
\end{theorem}

\begin{proof} 
First of all, it suffices to 
look at the $D_{f,g}$-cohomology of $V={\rm Fock}_{K\oplus N}$. 
Indeed, the proof of Proposition 8.2 of \cite{borvert} 
for the \emph{partial} lattice vertex 
algebra ${\rm Fock}_{M\oplus N}^\Sigma$
transfers easily to our case.

Now consider $v\in {\rm Ker}(D_{f,g})$ such that $H_Av=\alpha v$ 
for $\alpha<0$. Differentials $D_f$ and $D_g$ give $V$ the structure 
of the double complex, where the double grading
is provided by $(\deg\cdot,\cdot\deg^*)$. Differentials $D_f$ and $D_g$ 
increase this grading by $(0,1)$ and $(1,0)$ respectively.
Let $v=\sum_{a,b} v_{a,b}$ be the double-graded decomposition of $v$. 
It is enough to consider the case of $a+b=k$ for a fixed $k$, 
so we have $v=\sum_{i} v_i$ where the double grading of $v_i$ is 
$(k-i,i)$. The condition $D_{f,g}v=0$ implies that $D_g v_{i+1}=D_fv_{i}$.
Let $j$ be the maximum $i$ such that $v_i\neq 0$. We have $D_f v_j=0$,
so by Proposition \ref{dfprop}, there exists a $w$ such that $D_fw=v_j$.
Then $v-D_{f,g}w$ is again in the kernel of $D_{f,g}$, but it now has a smaller
value of $j$. By applying this procedure sufficiently many times,
we can get $j<0$, which implies that all $v_i$ are zero, since 
$\cdot\deg^*$ is nonnegative on $V$.

The $H_A=0$ statement is proved similarly. Let $D_{f,g}v=H_Bv=0$,
where $v=\sum_i v_i$ is as above. Let $j$ be the maximum value of 
$i$ such that $v_i\not\in {\rm Fock}_{K\oplus(K^*-\deg^*)}$.
Then we have $D_fv_j=D_g v_{j+1}\in {\rm Fock}_{K\oplus(K^*-\deg^*)}$.
Since $D_f$ preserves the $N$-grading, the $N$-graded part $v_j'$ of $v_j$
that is supported outside of $K^*-\deg^*$ satisfies $D_fv_j'=0$.
Then by Proposition \ref{dfprop} there exists $w\in V$ and such that 
$v_j'=D_fw$. The difference $v-D_{f,g}w$ will have smaller $j$, 
which eventually
leads to the situation where all 
$v_i\in {\rm Fock}_{K\oplus(K^*-\deg^*)}$.

It is clear from Corollary \ref{hb} that the eigenvalues of $H_A$
are integers.
Finally, the statements for $H_B$ 
are obtained by switching $K$ and
$K^*$.
\end{proof}

\begin{remark}\label{rfg}
Similarly, we can use Proposition \ref{somebigr} to show that 
for each $k$ $H_B=k$ eigenspace of $V_{f,g}$ comes from 
${\rm Fock}_{K\oplus (K^*-r\deg^*)}$ for a sufficiently big $r$.
\end{remark}

\begin{theorem}\label{sigmatype}
For strongly non-degenerate $f$ and $g$ the $N=2$ vertex algebra 
$V_{f,g}$ is of $\sigma$-model type.
\end{theorem}

\begin{proof}
It is clear that $H_A$ and $H_B$ are diagonalizable on $V_{f,g}$.
By Theorem \ref{main} the eigenvalues of $H_A$ and $H_B$
are nonnegative integers.
It has been shown in \cite[Lemma 4.5]{borlibg}
that the common eigenspaces of $L_{(1)}$ and $J_{(0)}$ of $V_{f,g}$ 
are finite-dimensional. Since $L_{(1)}\pm \frac12 J_{(0)}\geq 0$,
there are only finitely many such common eigenspaces for a fixed value 
of $L_{(1)}$, which proves that the eigenspaces of $L_{(1)}$ 
are finite-dimensional. It is also clear from Proposition \ref{lj}
that the eigenvalues of $L_{(1)}$ are in $\frac 12\ZZ$, which finishes 
the proof.
\end{proof}

\section{Relation to the conjectural description of
stringy cohomology}\label{sec.sh}
In this section we 
examine the structure of the chiral ring of $V_{f,g}$ in more detail
and connect it to the description of \emph{quantum stringy cohomology}
suggested in \cite{anvar}.

Let $K$ and $K^*$ be dual reflexive Gorenstein cones and let 
$f$ and $g$ be strongly non-degenerate coefficient functions.
Consider the ideal in the semigroup $\CC[K\oplus K^*]$ generated by 
$[m\oplus n], m\cdot n>0$. Denote the quotient by this ideal by
$\CC[L]$. Consider the space $V = \wedge(N_\CC)\otimes \CC[L]$
where, as usual, $\wedge$ means the exterior algebra.
We recall the following lemma.

\begin{lemma}\cite{anvar}
The space $V$ is equipped with a differential $d$ given by
$$
d:= \sum_{m} f(m) \contr m \otimes [m] +
\sum_{n} g(n) (\wedge n) \otimes [n]
$$
where $[m]$ and $[n]$ means multiplication by the corresponding monomials
in $\CC[K]\otimes\CC[K^*]$ acting on the module $\CC[L]$, $\wedge n$
means multiplication by $n$ in the exterior algebra and $\contr$
means contraction in the exterior algebra.
\end{lemma}\label{d}

\begin{remark}
In view of the isomorphism $\wedge(N_\CC)\cong \wedge(M_\CC)$
one could switch the roles of $M$ and $N$ in the construction of 
$V$ and $d$ without altering the resulting cohomology vector
space. This is related to \emph{spectral flow} isomorphism between
the $A$- and $B$-rings of $V_{f,g}$, see Remark \ref{remspectral}.
\end{remark}

It has been shown in \cite[Section 10]{anvar} that the cohomology of $V$
with respect to $d$ is isomorphic as a vector space to the conjectural
description of (quantum) stringy cohomology, in the case of a hypersurface
or a complete intersection induced by a nef-partition. We will now observe
that this space is also isomorphic to the $B$-ring of 
$V_{f,g}$. 

\begin{proposition}\label{secondmain}
Cohomology $H$ of $V$ with respect to $d$ is naturally isomorphic
to the $B$-ring of $V_{f,g}$. 
\end{proposition}

\begin{proof}
By Theorem \ref{main} and the definition of the $B$-ring
of $V_{f,g}$, it comes from  
${\rm Fock}_{(K-\deg)\oplus K^*}$. Since $D_{f,g}$ commutes with $H_B$,
the $B$-ring of $V_{f,g}$ is equal to the $D_{f,g}$ cohomology 
of the $H_B=0$ eigenspace of ${\rm Fock}_{(K-\deg)\oplus K^*}$.

By Proposition \ref{lj} and definition of $H_B$,
$$
H_Bv = ((\hat m+\deg)\cdot \hat n
+\sum_pi_p+\sum_qj_q+\sum_r(i_r+\frac12)+\sum_s(j_s-\frac12))v.
$$
As a result, the $H_B=0$ eigenspace of  ${\rm Fock}_{(K-\deg)\oplus K^*}$
is spanned by products of $n^{ferm}_{-\frac 12}$ and 
$\vert \hat m-\deg,\hat n\rangle$ with 
$\hat m\cdot \hat n=0$. So this eigenspace is isomorphic to $V$,
and it remains to show that the action of $d$ in Lemma \ref{d} 
is precisely the action of $D_{f,g}$. 

Let us calculate the action of 
$\oint m^{ferm}(z)\ee^{\int m^{bos}(z)}\,dz$
on $w\vert \hat m-\deg,\hat n\rangle$ where $w$ is a product of 
$n^{ferm}_{-\frac 12}$. 
Since the image will again satisfy $H_B=0$,
it has to be of the same form. As a result, the only mode of $m^{ferm}$
that can give a nonzero contribution will be $m^{ferm}_{\frac 12}$
that acts by a contraction on $w$. Similarly, we must have 
$m\cdot \hat n=0$, since otherwise we move out of $\Fock_{L}$.
Since we are taking $\oint$, the count of the degree of $z$
forces the corresponding mode from $\ee^{\int m^{bos}(z)}$
to be a shift $\gamma_m$.

The action of $\oint n^{ferm}(z)\ee^{\int n^{bos}(z)}\,dz$
on $w\vert \hat m-\deg,\hat n\rangle$ is calculated similarly. We need
to have $\hat m \cdot n=0$, and the mode of $n^{ferm}$ has to
be $n^{ferm}_{-\frac 12}$, which acts by multiplication. Then the 
count of the degree of $z$ forces us to have $\gamma_n$ for the 
mode of $\ee^{\int n^{bos}(z)}$. 

We have thus shown that $D_f+D_g$ and $d$ have identical actions on $V$,
which finishes the proof.
\end{proof}

In fact, by the results of \cite{anvar}, we can
quite explicitly calculate chiral rings in terms of the spaces 
$R_1(\theta,f)$ which we will now define.
Let $f$ and $g$ be non-degenerate coefficient functions.
For every face $\theta$ of the Gorenstein cone $K$ we will
abuse notations and denote the restriction of $f$ to $\theta$
by the same letter $f$. Then $f$ is non-degenerate for 
$\theta$ as well. We consider the ideal $I_f$ of $\CC[\theta]$
generated by the restrictions of the elements $z_i$ to $\CC[\theta]$.
Let $\CC[\theta^\circ]$ be the ideal of $\CC[\theta]$ which
corresponds to the interior of $\theta$. 
\begin{definition}\label{R1}
The natural inclusion $\CC[\theta^\circ]\to\CC[\theta]$
induces a (non-injective) map 
$$\CC[\theta^\circ]/I_f\CC[\theta^\circ]\to\CC[\theta]/I_f\CC[\theta]$$
and we denote its image by $R_1(\theta,f)$.
\end{definition}

\begin{remark}
Spaces $R_1(\theta,f)$ have been considered very early 
in the context of mirror symmetry, see \cite{Bat.var}.
The graded dimension of $R_1(\theta,f)$ has been calculated in
\cite{anvar}. In the case of simplicial $\theta$ it can be obtained
by a simple inclusion-exclusion formula in terms of
various $S$-polynomials of the faces of $\theta$. More generally, 
it also involves (intersection cohomology) $G$-polynomials of 
the partially ordered set 
of faces of $\theta$.
\end{remark}

For a face $\theta\subseteq K$ we denote by $\theta^*$ the
dual face of $K^*$, defined as the set of elements of $K^*$
which are zero on $\theta$.
\begin{theorem} \label{main.chiral}
Let $f$ and $g$ be strongly non-degenerate.
Then both chiral rings of the 
$N=2$ vertex algebra 
$V_{f,g}$ are naturally isomorphic {\bf as vector spaces} to
the space
$$
W_{f,g}=\oplus_{\theta\subseteq K} R_1(\theta,f) \otimes_\CC
R_1(\theta^*, g).
$$
\end{theorem}

\begin{proof}
We combine \cite[Theorem 10.2]{anvar} with Proposition \ref{secondmain}
for the $B$-ring case. The $A$-ring case differs by a switch 
of $M$ and $N$ and thus follows from the $B$-ring case.
\end{proof}

\begin{remark}\label{remchiral}
In fact, the isomorphism between the $A$- or $B$-ring and $W_{f,g}$
is pretty well-understood. We will state it for the $A$-ring.
Let $\Phi_\theta(z)$ be the 
field 
$(m_1)^{ferm}(z)\ldots (m_{\dim \theta})
^{ferm}(z)$ of ${\rm Fock}_{M\oplus N}$ 
where $m_1,\ldots,m_{\dim\theta}$ is a basis
of $\theta_\CC$. This is well-defined up to a constant
factor.
Then the part of the $A$-ring 
that corresponds to $\theta$ in $W_{f,g}$ 
comes from linear combinations of the fields of the 
form
$$
\Phi_\theta(z) \ee^{\int m^{bos}(z)+n^{bos}(z)-(\deg^*)^{bos}(z)}
$$
for $m\in\theta^\circ$ and $n\in(\theta^*)^\circ$.
The proof of this statement follows from the argument
of \cite[Theorem 10.2]{anvar}. It is rather non-trivial,
as we will see from the example below.
\end{remark}

\begin{example}\label{example.quintic}
Let $\Delta_1^*$ be the reflexive polytope in $\ZZ^4$ with vertices
$\hat n_1=(1,0,0,0)$, $\hat n_2=(0,1,0,0)$, $\hat n_3=(0,0,1,0)$, 
$\hat n_4=(0,0,0,1)$, 
and $\hat n_5=(-1,-1,-1,-1)$. The only other lattice point in $\Delta_1^*$
is $\hat n_0=(0,0,0,0)$. We remark that this corresponds to the 
famous quintic example. 
We consider the corresponding 
cone $K^*\subseteq \ZZ^4\oplus \ZZ$ and introduce $n_i=(\hat n_i,0)$.
The convex hull of $\{n_i\}$ will be denoted by $\Delta^*$.
Notice that $n_0=\deg^*$.
It is rather straightforward to see that in this 
case $R_1(\theta^*,g)$ is zero unless $\theta^*=\{0\}$ or 
$\theta^*=K^*$.
The space $R_1$ for zero-dimensional face is one-dimensional,
so according to Theorem \ref{main.chiral} 
the $A$- and $B$-rings are isomorphic as
vector spaces to $R_1(K,f)\oplus R_1(K^*,g)$.

We will only describe the $R_1(K^*,g)$ part inside the 
$A$-ring, which turns out to be a four-dimensional vector space. 
It corresponds to $\theta=\{0\}$, so according to 
Remark \ref{remchiral}, 
we need to look at linear combinations of the 
fields of ${\rm Fock}_{M\oplus N}$ 
of the form 
\begin{equation}\label{diag.fields}
\ee^{\int n^{bos}(z)-n_0^{bos}(z)}
\end{equation}
for $n\in (K^*)^\circ$. It is clear that all such fields 
commute with the differential $D_{f,g}$ and therefore descend 
to $V_{f,g}$. It is also clear that they satisfy $H_A=0$
and hence descend to elements of the $A$-ring.

Let us fix a non-degenerate coefficient function $g$, i.e.
six complex numbers $g(n_i),~i=0,\ldots,5$. 
For every $m\in M$ and every $n\in K^*$ the field
$m^{ferm}(z)\ee^{\int n^{bos}(z)}$ anticommutes with 
$D_f$ and gives 
$$
\sum_{i=0}^5 g(n_i)(m\cdot n_i)\ee^{\int n_i^{bos}(z)+n^{bos}(z)}
$$
when acted upon by $D_{f,g}$. As a result, if we identify
the space of fields of the form \eqref{diag.fields} with
$\CC[(K^*)^\circ]$, we see that the map to the $A$-ring
passes through the quotient $\CC[(K^*)^\circ]/I_{g}
\CC[(K^*)^\circ]$. In this particular example, this 
means that it is enough to consider the linear combinations
of the fields $\ee^{\int kn_0^{bos}(z)}$ for $k=0,1,2,3,4$.
We leave some easy commutative algebra calculations to 
the reader.

It is much more difficult to see that the map from 
$\CC[(K^*)^\circ]$ to the $A$-ring 
passes through $R_1(K^*,g)$. This means that in fact
$\ee^{\int 4n_0^{bos}(z)}$ maps to zero in the $A$-ring.
Denote by $m_i,~i=1,\ldots,5$ the elements of the cone $K$ 
that correspond to vertices of the dual polytope $\Delta$. 
More specifically,
they will be given by $(4,-1,-1,-1,1)$, $(-1,4,-1,-1,1)$,
$(-1,-1,4,-1,1)$, $(-1,-1,-1,4,1)$ and $(-1,-1,-1,-1,1)$ respectively.
Their scalar products with $n_j$ are $(m_i\cdot n_j)=5\delta_i^j$
for $j>0$ and $(m_i\cdot n_0)=1$.

Consider the fields
$$
\begin{array}{l}
R_1(z) =
\frac {1}{5g_1} m_1^{ferm}(z)
\ee^{\int (n_2+n_3+n_4+n_5-n_0)^{bos}(z)}
\\
R_2(z)=
\frac {g_0}{5^2g_1g_2} m_2^{ferm}(z)
\ee^{\int (n_3+n_4+n_5)^{bos}(z)}
\\
R_3(z)=
\frac {g_0^2}{5^3g_1g_2g_3} m_3^{ferm}(z)
\ee^{\int (n_4+n_5+n_0)^{bos}(z)}
\\
R_4(z)=
\frac {g_0^3}{5^4g_1g_2g_3g_4} 
m_4^{ferm}(z)
\ee^{\int (n_5+2n_0)^{bos}(z)}
\\
R_5(z)=
\frac {g_0^4}{5^5g_1g_2g_3g_4g_5} 
m_5^{ferm}(z)
\ee^{\int 3n_0^{bos}(z)}.
\end{array}
$$
We observe that $D_f R_i(z)= 0$ for all $i=1,\ldots, 5$.
Indeed, when we write OPE of $m^{ferm}(z) \ee^{\int m^{bos}(z)}$
with $R_1(w)$ for some $m\in \Delta$, we will only get a pole
at $z=w$ for the bosonic part when $m=m_1$. However, in this case
the fermions contribute a zero of order one which cancels 
the singularity.  There are no poles in the bosonic parts 
of OPE for any of the other $R_i$, since the exponents lie
in $K^*$.

On the other hand, when we calculate $D_g R_i(z)$ only 
two terms of $D_g$ will matter, namely the terms for $n_i$ and 
$n_0$, in view of the scalar products $m_i\cdot n_j$. This yields
$$D_g R_i(z) = 
\frac{g_0^{i-1}}{5^{i-1}g_1\cdots g_{i-1}}
\ee^{\int (n_{i}+\ldots + n_5+(i-2)n_0)^{bos}(z)} 
$$$$+
\frac{g_0^{i}}{5^ig_1\cdots g_i}
\ee^{\int (n_{i+1}+\ldots + n_5+(i-1)n_0)^{bos}(z)}.
$$
As a result, $D_{f,g}\sum_{i=1}^{5} (-1)^{i-1} R_i(z)$ 
gives $(1 + \frac {g_0^5}{5^5g_1\cdot g_5})\ee^{\int 4n_0^{bos}(z)}$
where we have used the fact that $\sum_{i=1}^5 n_i = 5 n_0$.
This shows that $\ee^{\int 4n_0^{bos}(z)}$ maps to zero in $V_{f,g}$.
It is amusing to notice that the coefficient in front is zero exactly
when the dual to the quintic has additional singularities, i.e.
when $g$ is degenerate.

As a result of the above calculation we see that the diagonal part 
of the quantum cohomology ring of the quintic is isomorphic to 
$\CC[t]/t^4$, where $t$ is the image of $\ee^{\int n_0^{bos}(z)}$.
\end{example}

\begin{remark}\label{remspectral}
We can observe that the part of $W_{f,g}$ that correspond 
to $\theta$ comes from 
$$
\Phi_\theta(z) \ee^{\int m^{bos}(z)+n^{bos}(z)-(\deg^*)^{bos}(z)},
$$
for the $A$-ring and comes from
$$
\Psi_{\theta^*}(z) \ee^{\int m^{bos}(z)+n^{bos}(z)-\deg^{bos}(z)},
$$
for the $B$-ring, where $\Psi$ is defined analogously to
$\Phi$. This amounts to an action of the field
$$S(z)=\Phi_K(z)\ee^{\int \deg^{bos}(z) - (\deg^*)^{bos}(z)}.$$
We can compare this with the spectral flow defined physically
in \cite{LVW}, by noticing that \emph{formally}
$S(z)=\ee^{\int J(z)}$, in view of the boson-fermion correspondence
for the fermionic part of ${\rm Fock}_{M\oplus N}$.
\end{remark}

To end this section, we remark on the relation to
non-quantum stringy cohomology spaces
of the Calabi-Yau hypersurfaces in toric varieties.
These have been defined in \cite{anvar} as follows.
Let $\Delta$ and $\Delta^*$ be dual reflexive polytopes 
in $M_1$ and $N_1=M_1^*$ and let
$K$ and $K^*$ be the corresponding reflexive Gorenstein cones 
in the lattices $M=M_1\oplus \ZZ\deg$ 
and $N=N_1\oplus \ZZ\deg^*$. Let $\Sigma_1$
be a regular fan in $N_1$ whose one-dimensional cones contain all
vertices of $\Delta^*$. We can define a decomposition $\Sigma$ of $N$ into
a union of semigroups by extending $\Sigma'$ into the $\ZZ\deg^*$ 
direction. When intersected with $K^*$, $\Sigma$ becomes a fan,
which we will denote by the same letter. Consider the \emph{partial
semigroup ring} $\CC[M\oplus N]^\Sigma$. It is isomorphic as 
a vector space to the ring $\CC[M\oplus N]$ but the multiplication product
is redefined by 
$$
[m\oplus n][\hat m\oplus \hat n]=\Big\{
\begin{array}{ll}
[(m+\hat m)\oplus(n+\hat n)],
&{\rm if~} n,\hat n\in C,{\rm~for~some~}C\in\Sigma\\
0,&{\rm otherwise.}
\end{array}
$$
We redefine the shifts $\gamma_{m\oplus n}$ as multiplication 
in the partial semigroup ring. 
Then the stringy cohomology for the given two
coefficient functions is again the cohomology of the differential 
of Lemma \ref{d}. 

One can use the redefined shifts $\gamma_{m\oplus n}$ to
construct \emph{partial lattice vertex algebra} 
${\rm Fock}_{M\oplus N}^\Sigma$ (see \cite{borvert}).
Then one can still define coefficient functions and the differential
$D_{f,g}$.
We denote the cohomology of ${\rm Fock}_{M\oplus N}^\Sigma$
by $D_{f,g}$ by $V_{f,g}^\Sigma$. 
We claim that its chiral
rings again equal the stringy cohomology. Indeed, the calculations
for the $A$-rings are unchanged, since the differential
$D_f$ is not affected by $\Sigma$, and Proposition \ref{dfprop}
still holds. It is slightly more difficult to show that the $B$-ring 
of $V_{f,g}^\Sigma$ can be calculated by the differential 
$D_{f,g}$ on the $H_B=0$ part of 
${\rm Fock}_{(K-\deg)\oplus K^*}^{\Sigma}$. 
Indeed, we can not
simply follow the proof of  Proposition \ref{dfprop} with the 
roles of $K$ and $K^*$ interchanged, because we need to be able 
to degenerate $\CC[K^*]^\Sigma$ according to a fan $\Phi$, which
is impossible. 
However, one can show that $V_{f,g}^\Sigma$ is isomorphic to
$D_{f,g}$ cohomology of the space ${\rm Fock}_{(K-\deg)\oplus N}$
by modifying the argument of  \cite[Proposition 8.2]{borvert}
slightly. Then one can prove that $H_B$ is nonnegative
on $V_{f,g}$ and $(H_B=0)$-part of $D_f$-cohomology 
of ${\rm Fock}_{(K-\deg)\oplus N}^{\Sigma}$ comes from
${\rm Fock}_{(K-\deg)\oplus K^*}^{\Sigma}$, by following the 
arguments of Proposition \ref{dfprop}. More specifically,
we now assume that $n_0\not\in K^*$, which means that there
exists a vertex $\hat m$ of $\Delta$ with $\hat m\cdot n < 0$.
The formula \eqref{lm} for $l_m$ is unchanged. We still basically
have a Koszul complex, and can perform the degeneration
argument with respect to a triangulation $\Phi$. Instead of 
${\rm Box}(C)$ we will be now dealing with ${\rm Box}(C)-\deg$.
The inequality
\eqref{technical} will then translate into
$$
H_B=(\sum_ia_im_i)\cdot n_0 - \sum_{l_i<0} (l_i-\frac 12)
$$
$$
=(\sum_ia_im_i)\cdot n_0 - \sum_{l_i<0} m_i\cdot n_0
\geq 
\sum_{l_i<0} (a_i-1)m_i\cdot n_0
$$
which will be positive due to the term for $m_i=\hat m$.
The transition
from $D_f$-cohomology to $D_{f,g}$ cohomology in Theorem \ref{main}
is
unchanged for the partial semigroup case. 

\begin{remark}
The dimensions of the chiral rings are clearly
unchanged after passing to $V_{f,g}^\Sigma$, in view of the 
calculation of \cite{anvar} that shows that the dimensions
of $R_1$ are unchanged.
\end{remark}

\begin{remark}
We have not mentioned the double grading on $A$-rings and $B$-rings 
of $V_{f,g}$ and $V_{f,g}^\Sigma$, which corresponds to
the Hodge structure $H=\oplus H^{p,q}$ 
on the stringy cohomology. The Hodge components 
are the common eigenspaces of $J_{(0)}$ and the grading
$\deg\cdot+\cdot \deg^*$ which clearly descends to the cohomology
of $D_{f,g}$. It appears that in some vague sense, 
this grading is what's left of anti-holomorphic fields 
of the vertex algebra in the definition of Kapustin and Orlov
\cite{KO} after taking $\bar \partial$ cohomology.
However we do not know how to construct a Kapustin-Orlov vertex
algebra responsible for mirror symmetry.
\end{remark}

\section{Open questions}\label{sec.oq}
Hopefully, the reader is convinced that the $N=2$ vertex algebras
$V_{f,g}$ provide rich examples of vertex algebras.
However, many of their properties are still not well-understood.
We list here several open questions of algebraic nature
related to this construction,
but this list is by no means complete.

\begin{question}
Are $V_{f,g}$ actually different as \emph{vertex algebras}?
This is important, since a negative answer would perhaps provide 
a non-trivial connection on the chiral rings of $V_{f,g}$. Of course,
one has to look at generic pairs $(f,g)$.
\end{question}

\begin{question}
Is there a way to see GKZ system of differential equations \cite{GKZ} 
in the context of $V_{f,g}$? This question seems related to the 
question of a connection on the space of $V_{f,g}$-s considered
as a bundle over the space of the parameters $(f,g)$.
\end{question}

\begin{question}
Are algebras $V_{f,g}$ generated by a finite number of fields?
While some examples suggest this, we are not yet able to prove it
in general.
\end{question}

\begin{question}
It is not hard to show that for a given $k$ there is a Zariski open
set in the space of coefficient functions $(f,g)$ such that 
the $L_{(1)}$ eigenspace of $V_{f,g}$ has a constant dimension.
Is it true that there is a Zariski open set that works for all $k$
simultaneously? This is the question of whether $V_{f,g}$ form
a flat family of vertex algebras over an algebraic parameter space.
\end{question}

\begin{question}
The results of this paper provide the quantum and non-quantum
stringy cohomology spaces of \cite{anvar} 
with a ring structure. Is it possible to describe this structure
without a reference to vertex algebras? Note that it is 
quite easy to see this product structure for the ``diagonal part''
as in Example \ref{example.quintic}, but the general case
is at the moment open.
\end{question}

\end{document}